# The theory of opetopes via Kelly-Mac Lane graphs

Eugenia Cheng

Department of Pure Mathematics, University of Cambridge
E-mail: e.cheng@dpmms.cam.ac.uk

October 2002

**Abstract**

This paper follows from the work of [4] and [3]. In the former paper we give an explicit construction of opetopes, the underlying cell shapes in the theory of opetopic $n$-categories; at the heart of this construction is the use of certain trees. In the latter paper we give a description of trees using Kelly-Mac Lane graphs. In the present paper we apply the latter to the former, to give a construction of opetopes using Kelly-Mac Lane graphs.

## Contents

## Introduction

In [1] Baez and Dolan give a definition of weak $n$-category in which the underlying shapes of cells are 'opetopes' and the underlying data are given by 'opetopic sets'. In the present paper we give an alternative construction of opetopes using the language of closed categories. This is made possible by the results of [3], giving a precise correspondence between the trees involved in the construction of opetopes, and the allowable morphisms in certain closed categories.

The end result is a precise algebraic characterisation in place of the more geometric world of trees.

The idea is as follows. Recall ([4]) that a $k$-opetope is constructed as a 'configuration' for composing $(k-1)$-opetopes; this is expressed as a tree (see [4]) whose nodes are labelled by the $(k-1)$-opetopes in question, with the edges giving their inputs and outputs.

In order to express this more precisely, it is helpful to give a more formal description of trees, as a bijection between inputs and outputs of nodes subject to a condition ensuring that closed loops do not arise.

In fact, this leads to an abstract description of trees as certain Kelly-Mac Lane graphs; this is the subject of [3]. In the present paper we apply the results of [3] to the construction of opetopes.

First we need to generalise the work of [3] since we require the 'labelled' version of the trees presented in that paper, and this is the subject of Section 1.1.

In Section 1.2 we give the actual construction. We use the result, from [3], that a tree is precisely an allowable Kelly-Mac Lane graph of a certain shape. We may therefore express the process of forming labelled trees more precisely, by seeking allowable morphisms in a certain closed category. The construction we use to build up dimensions in this way is what we call a 'ladder'; the idea is to pick out precisely the allowable morphisms that satisfy two conditions. The first is that the morphism must be of the correct shape to correspond to a tree. For the second, recall that an arrow of a slice multicategory is a way of composing its source elements to give its target element; the second condition for the 'ladder' construction corresponds to this stipulation.

In Section 1.3 we give some low-dimensional examples of the above construction to help elucidate this rather abstract approach.

Then, in Section 2, we prove that this construction does give the same opetopes as defined in [4]. This is the main result of this work. We conclude, in Section 3, with a brief discussion about the category of opetopes. This category is constructed explicitly in [2], and the work in the present paper extends to the construction of this category. We do not explicitly give this construction here, but we give some low-dimensional examples of the 'face' maps of opetopes in the new framework.

**Terminology**

i) Since we are concerned chiefly with *weak* $n$-categories, we follow Baez and Dolan ([1]) and omit the word 'weak' unless emphasis is required; we refer to *strict* $n$-categories as 'strict $n$-categories'.

ii) In [1] Baez and Dolan use the terms 'operad' and 'types' where we use 'multicategory' and 'objects'; the latter terminology is more consistent with Leinster's use of 'operad' to describe a multicategory whose 'objects-object' is 1.

iii) As in [2], we work with opetopes not precisely the same as those given in [1], but rather as given by the opetopic theory as explained in [4] using multicategories with a category rather than a set of objects; we refer the reader to [4] for the full details.

**Acknowledgements**

This work was supported by a PhD grant from EPSRC. I would like to thank Martin Hyland and Tom Leinster for their support and guidance.

# 1 Opetopes

In this section we use the results of [3] to construct opetopes. We first need to introduce the notion of *labelled* Kelly-Mac Lane graphs.

## 1.1 Preliminaries

For the construction of opetopes we require the 'labelled' version of the theory presented in [3]: labelled shapes, labelled graphs and labelled trees.

Given a category $\mathbb{C}$ we can form *labelled shapes* (in $\mathbb{C}$), that is, shapes labelled by the objects of $\mathbb{C}$. A labelled shape is a shape $T$ with each 1 'labelled' by an object $A_i$ of $\mathbb{C}$. We write this as

$$|T|(A_1, \ldots, A_k).$$

The variable set is then defined to be the variable set of the underlying shape.

For example given

$$T = [\ [1,1] \otimes 1 \otimes 1 \ , \ I \ ] \otimes 1$$

we have a labelled shape

$$\alpha = |T|(A_1, \ldots, A_5) = [\ [A_1, A_2] \otimes A_3 \otimes A_4 \ , \ I \ ] \otimes A_5$$

with underlying shape $T$, and

$$v(\alpha) = v(T) = \{+, -, -, -, +\}.$$

Given a category $\mathbb{C}$ we can form labelled graphs, that is, graphs whose edges are labelled by morphisms of $\mathbb{C}$ as follows. Consider labelled shapes $\alpha$ and $\beta$ with underlying shapes $T$ and $S$ respectively. A *labelled graph*

$$\alpha \longrightarrow \beta$$

is a graph
$$\xi : T \longrightarrow S$$
together with a morphism $x \longrightarrow y$ for each pair of labels $x, y$ whose underlying variables are mates under $\xi$, with $v(x) = -$ and $v(y) = +$ in the twisted sum. That is, the morphism is in the direction
$$- \longrightarrow +.$$
For example, the following is a labelled graph

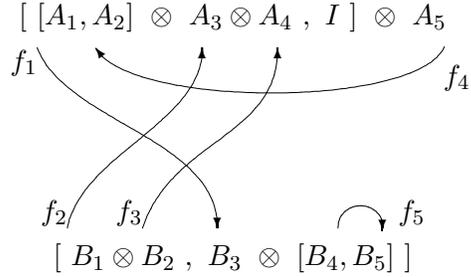

with underlying graph and variances as below

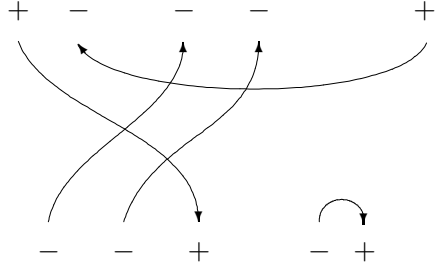

.

Observe that, since the variances of the domain are reversed in the twisted sum, the direction of morphisms is also reversed in the domain.

We write $K\mathbb{C}$ for the category of labelled shapes and labelled graphs in $\mathbb{C}$; thus $G = K\mathbf{1}$ as mentioned in [3].

A labelled graph is called allowable if and only if its underlying graph is allowable. We write $A\mathbb{C}$ for the category of labelled shapes and allowable labelled morphisms. We observe immediately that the correspondence between trees and graphs exhibited in [3] generalises to a correspondence between labelled graphs and labelled trees.

**Lemma 1.1** *A labelled tree in $\mathbb{C}$ is precisely an allowable morphism*
$$\alpha_1, \otimes \cdots \otimes \alpha_k \longrightarrow \alpha \in A\mathbb{C}$$

*with underlying shape*

$$X_{m_1} \otimes \cdots \otimes X_{m_k} \longrightarrow X_{(\sum_i m_i)-k+1}.$$

Recall ([4]) that a labelled tree gives a 'configuration for composing' arrows of a symmetric multicategory via object-morphisms, as used in the slice construction. By the above correspondence, a labelled graph as above may also be considered to give such a configuration; thus in Section 2.1 we will use such graphs to give an alternative description of the slice construction. We will need the following construction.

Given categories $\mathbb{C}$ and $\mathbb{D}$ and a functor

$$F : \mathbb{C} \longrightarrow K\mathbb{D}$$

we define a functor

$$KF : K\mathbb{C} \longrightarrow K\mathbb{D}$$

as follows.

- on objects

An object in $K\mathbb{C}$ is a labelled shape

$$\alpha = |T|(x_1, \ldots, x_n);$$

put

$$KF(\alpha) = |T|(Fx_1, \ldots Fx_n) \in K\mathbb{D}.$$

- on morphisms

Given a morphism

$$|T|(x_1, \ldots, x_n) \xrightarrow{f} |S|(x_{n+1}, \ldots, x_m) \in K\mathbb{C}$$

we define $KFf$ as follows. Suppose $f$ has underlying graph $\xi$, say. Consider a pair of mates $a$ and $b$ in $\xi$, with the edge joining them 'labelled' with morphism

$$g : a \longrightarrow b \in \mathbb{C}.$$

This gives a morphism

$$Fg : Fa \longrightarrow Fb \in K\mathbb{D}.$$

So $Fg$ is a graph labelled in $\mathbb{D}$. Then $KFf$ consists of all such graphs given by mates in $\xi$, positioned according to the positions in $\xi$.

Furthermore, if $F : \mathbb{C} \longrightarrow \mathbb{D}$ then we get

$$AF : A\mathbb{C} \longrightarrow A\mathbb{D}$$

by restricting the functor $KF$.

## 1.2 The construction of opetopes

We seek to define, for each $k \geq 0$, a category $\mathbf{Ope}_k$ of $k$-opetopes. A $k$-opetope $\theta$ is to have a list of input $(k-1)$-opetopes $\alpha_1, \ldots, \alpha_m$, say, and an output $(k-1)$-opetope $\alpha$, say. This data is to be expressed as an object

$$[\,\alpha_1 \otimes \cdots \otimes \alpha_m\,,\,\alpha\,] \in A\mathbf{Ope}_{k-1}$$

called the *frame* of $\theta$ (see [1]). Each frame has shape $X_m = [1^{\otimes m}, 1]$ for some $m \geq 0$. So, for each $k$ we will have a functor

$$\phi_k : \mathbf{Ope}_k \longrightarrow A\mathbf{Ope}_{k-1}$$

and thus

$$A\phi_k : A\mathbf{Ope}_k \longrightarrow A\mathbf{Ope}_{k-1}.$$

$\mathbf{Ope}_k$ is defined inductively; for $k \geq 2$ it is a certain full subcategory of the comma category

$$(I \downarrow A\phi_{k-1})$$

with the following motivation. A $k$-opetope $\theta$ with frame

$$[\,\alpha_1 \otimes \cdots \otimes \alpha_m\,,\,\alpha\,]$$

is a configuration for composing $\alpha_1, \ldots, \alpha_m$ to result in $\alpha$. That is, it is an allowable morphism

$$I \xrightarrow{\theta} [\phi_{k-1}\alpha_1 \otimes \cdots \otimes \phi_{k-1}\alpha_m\,,\,\phi_{k-1}\alpha\,] \in A\mathbf{Ope}_{k-2}$$

such that the composition does result in $\alpha$. Such a $\theta$ is clearly an object of $(I \downarrow A\phi_{k-1})$; so we take the full subcategory whose objects are all those of the correct form.

In fact we begin with a more general construction for building up dimensions.

**Definition 1.2** *A* ladder *is given by*

- *for each $k \geq 0$ a category $\mathbb{D}_k$*
- *for each $k \geq 1$ a functor $F_k : \mathbb{D}_k \longrightarrow A\mathbb{D}_{k-1}$*

*such that for each $k \geq 2$, $F_k$ is of the form*

$$\mathbb{D}_k \longrightarrow (I \downarrow AF_{k-1}) \longrightarrow A\mathbb{D}_{k-1}$$

*where the second morphism is the forgetful functor.*

Note that given $F_k : \mathbb{D}_k \longrightarrow A\mathbb{D}_{k-1}$ we have a functor

$$AF_k : A\mathbb{D}_k \longrightarrow A\mathbb{D}_{k-1}$$

and the comma category $(I \downarrow AF_{k-1})$ has as its objects pairs $(\theta, z)$ where $z \in A\mathbb{D}_{k-1}$ and $\theta$ is an allowable morphism

$$I \xrightarrow{\theta} AF_{k-1}(z) \in A\mathbb{D}_{k-1}.$$

**Definition 1.3** *The* opetope ladder *is given as follows.*

- $\mathbb{D}_0 = 1 = \{x\}$, say

- $\mathbb{D} = 1 = \{u\}$, say, with

$$\phi_1 : \mathbb{D}_1 \longrightarrow A\mathbb{D}_0$$
$$u \longmapsto [x, x]$$

- For $k \geq 2$, $\mathbb{D}_k$ is a full subcategory of $(I \downarrow A\phi_{k-1})$. This comma category has objects $(\theta, z)$ where $z \in A\mathbb{D}_{k-1}$ and

$$I \xrightarrow{\theta} A\phi_{k-1}(z)$$

  is an allowable morphism in $A\mathbb{D}_{k-2}$. Then the subcategory $\mathbb{D}_k$ by the following two conditions:

  A. The objects of $\mathbb{D}_k$ are all $(\theta, z)$ such that $z$ has shape $X_m$ for some $m \geq 0$. So $z = [\alpha_1 \otimes \cdots \otimes \alpha_m, \alpha]$ for some $\alpha_i, \alpha \in \mathbb{D}_{k-1}$.

  B. For $k \geq 3$ we require in addition that

  $$A\phi_{k-2}\bar{\theta} \circ (\alpha_1 \otimes \cdots \otimes \alpha_m) = \alpha$$

  as morphisms in $A\mathbb{D}_{k-3}$.

- For $k \geq 2$ the functor $\phi_k : \mathbb{D}_k \longrightarrow A\mathbb{D}_{k-1}$ is the following composite

$$\mathbb{D}_k \hookrightarrow (I \downarrow A\phi_{k-1}) \longrightarrow A\mathbb{D}_{k-1}$$

  where the functors shown are the inclusion followed by the forgetful functor.

Note that the composition in condition B is possible: each $\alpha_i$ is an object of $\mathbb{D}_{k-1}$, so is by definition a morphism

$$I \longrightarrow A\phi_{k-2}(\phi_{k-1}\alpha_i) \in A\mathbb{D}_{k-3}.$$

Now $\theta$ is a morphism

$$I \longrightarrow [\phi_{k-1}\alpha_1 \otimes \cdots \otimes \phi_{k-1}\alpha_m \,,\, \phi_{k-1}\alpha \,]$$

so

$$\bar{\theta} : \phi_{k-1}\alpha_1 \otimes \cdots \otimes \phi_{k-1}\alpha_m \longrightarrow \phi_{k-1}\alpha$$

so the domain of $A\phi_{k-2}\bar{\theta}$ is indeed the codomain of $(\alpha_1 \otimes \cdots \otimes \alpha_m)$ and the composite in $A\mathbb{D}_{k-3}$ may be formed.

**Definition 1.4** *For each $k \geq 0$ the category $\mathbb{D}_k$ defined above is the category of k-opetopes. We write $\mathbb{D}_k = \mathbf{Ope}_k$. If the frame of a k-opetope has shape $X_m$ we say $\theta$ is an m-ary opetope.*

**Remarks 1.5**

1) In general (that is for $k \geq 3$) the objects of $\mathbb{D}_k$ are those of $(I \downarrow A\phi_{k-1})$ satisfying the conditions A and B. Condition A restricts our scope only to those objects having the correct shape; condition B ensures that the 'output' of the opetope is indeed the composite given. For $k = 2$ condition B does not apply; any configuration of composing identity maps gives the identity.

2) A morphism $(\theta, z) \xrightarrow{f} (\theta', z')$ in $(I \downarrow A\phi_{k-1})$ is a morphism
$$f : z \longrightarrow z' \in A\mathbb{D}_{k-1}$$
such that the following diagram commutes:

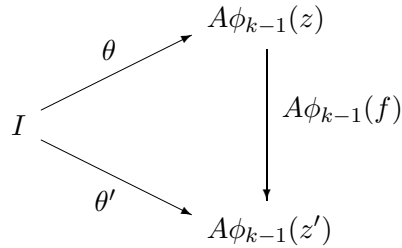

so a morphism $\theta \xrightarrow{f} \theta'$ in $\mathbb{D}_k$ is given as follows. Writing
$$\phi_k\theta = [\alpha_1 \otimes \cdots \otimes \alpha_m, \alpha]$$
$$\phi_k\theta' = [\beta_1 \otimes \cdots \otimes \beta_j, \beta]$$

$f$ must be a morphism
$$[\alpha_1 \otimes \cdots \otimes \alpha_m, \alpha] \longrightarrow [\beta_1 \otimes \cdots \otimes \beta_j, \beta] \in A\mathbb{D}_{k-1}.$$

So we must have $m = j$ and $f$ has the form

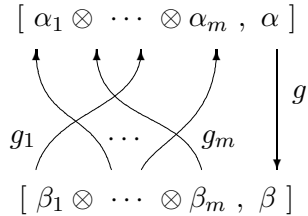

that is, a permutation $\sigma \in \mathbf{S}_m$ and morphisms

$$g_i : \beta_i \longrightarrow \alpha_{\sigma(i)}, \quad \text{for each } 1 \leq i \leq m$$

$$g : \alpha \longrightarrow \beta$$

in $\mathbb{D}_{k-1}$, such that the following diagram commutes

$$\begin{array}{c}
\xymatrix{
& [\ \phi_{k-1}\alpha_1 \otimes \cdots \otimes \phi_{k-1}\alpha_m\ ,\ \phi_{k-1}\alpha\ ] \ar[d]^{\phi_{k-1}g} \\
I \ar[ur]^{\theta} \ar[dr]_{\theta'} & \\
& [\ \phi_{k-1}\beta_1 \otimes \cdots \otimes \phi_{k-1}\beta_m\ ,\ \phi_{k-1}\beta\ ]
}
\end{array}$$

### 1.3 Examples

We now give the first few stages of the construction explicitly, together with some examples.

- $k = 0$

$\mathbf{Ope}_0 = \mathbf{1}$, that is, there is only one 0-opetope. We may think of this as an object $\cdot$ ; we write $x$ for convenience.

- $k = 1$

$\mathbf{Ope}_1 = \mathbf{1}$, that is, there is only one 1-opetope, $u$, say. We have

$$\phi_1(u) = [x, x] \in A\mathbf{Ope}_0$$

that is, the unique 1-opetope $u$ has one input 0-opetope and one output 0-opetope. We may think of this as

$$\longrightarrow$$

or, showing variances

$$-\ \ +$$

and then we have

$$\phi_1(1_u) \quad = \quad \begin{array}{c} - \\ | \\ | \\ - \end{array} \ \begin{array}{c} + \\ | \\ | \\ + \end{array} \quad ,$$

an allowable morphism in $A\mathbf{Ope}_0$. (We do not show arrowheads since all arrows in $\mathbf{Ope}_0$ are identity arrows.)

- $k = 2$

We seek to construct the category $\mathbf{Ope}_2$. First we consider an object $\alpha \in \mathbf{Ope}_2$. $\alpha$ has frame
$$\phi_2 \alpha \in A\mathbf{Ope}_1$$
where $\phi_2 \alpha$ has shape $X_m$ for some $m \geq 0$. So we have
$$\phi_2 \alpha = [u^{\otimes m}, u] = [u \otimes \cdots \otimes u, u].$$

Now $\alpha$ is an allowable morphism
$$I \xrightarrow{\alpha} [\, \phi_1 u \otimes \cdots \otimes \phi_1 u \,,\, \phi_1 u \,] \in A\mathbf{Ope}_0 = A\mathbf{1}$$
that is
$$I \xrightarrow{\alpha} [\, [x,x] \otimes \cdots \otimes [x,x] \,,\, [x,x] \,]$$
or equivalently a morphism
$$[x,x] \otimes \cdots \otimes [x,x] \longrightarrow [x,x] \in A\mathbf{1}.$$

For example, for $m = 3$ we may have a graph

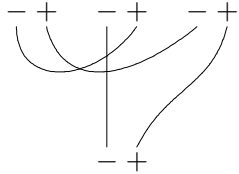

which we will later see corresponds to the following

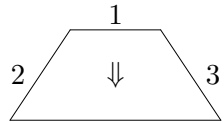

where the numbers show the order in which the input 1-opetopes are listed.

For the nullary case $m = 0$ we seek an allowable morphism
$$I \longrightarrow [x,x].$$
There is precisely one such, given by the following graph and we will later see that this corresponds to the nullary 2-opetope

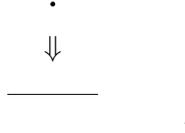

We now consider a morphism
$$\alpha \xrightarrow{f} \alpha' \in \mathbf{Ope}_2.$$
We must have
$$\phi_2\alpha = \phi_2\alpha' = [u^{\otimes m}, u],$$
say. Then $f$ is a morphism
$$[u^{\otimes m}, u] \longrightarrow [u^{\otimes m}, u] \in A\mathbf{Ope}_1 = A\mathbf{1}.$$
So $f$ must be a permutation $\sigma \in \mathbf{S}_m$, an isomorphism. So we have
$$\mathbf{Ope}_2(\alpha, \alpha') = \begin{cases} \mathbf{S}_m & \text{if } \alpha \text{ and } \alpha' \text{ are both } m\text{-ary} \\ \emptyset & \text{otherwise} \end{cases}$$
and $\mathbf{Ope}_2$ is equivalent to a discrete category whose objects are the natural numbers.

Note that the action of $\phi_2$ on morphisms is given as follows. Given a morphism $f$ as above, the morphism
$$\phi_2 f : \phi_2\alpha \longrightarrow \phi_2\alpha' \in A\mathbf{Ope}_1$$
is given by the forgetful functor
$$(I \downarrow A\phi_1) \longrightarrow A\mathbf{Ope}_1$$
so is simply the graph given by the permutation $\sigma$.

- $k = 3$

We now seek to construct the category $\mathbf{Ope}_3$. We first consider an $m$-ary opetope $\theta \in \mathbf{Ope}_3$ with frame
$$[\alpha_1 \otimes \cdots \otimes \alpha_m, \alpha] \in A\mathbf{Ope}_2$$

such that
$$\phi_2\alpha_i = [u^{\otimes n_i}, u] \quad \text{for each } 1 \leq i \leq m$$
$$\phi_2\alpha = [u^{\otimes n}, u].$$

So $\theta$ is an allowable morphism
$$I \xrightarrow{\theta} [\, [u^{\otimes n_1}, u] \otimes \cdots \otimes [u^{\otimes n_m}, u], [u^{\otimes n}, u]\,]$$

or equivalently
$$[u^{\otimes n_1}, u] \otimes \cdots \otimes [u^{\otimes n_m}, u] \xrightarrow{\bar{\theta}} [u^{\otimes n}, u] \in A\mathbf{Ope}_1,$$

such that
$$(A\phi_1)\bar{\theta} \circ (\alpha_1 \otimes \cdots \otimes \alpha_m) = \alpha$$

as morphisms in $A\mathbf{Ope}_0$.

For example for $m = 2$ consider

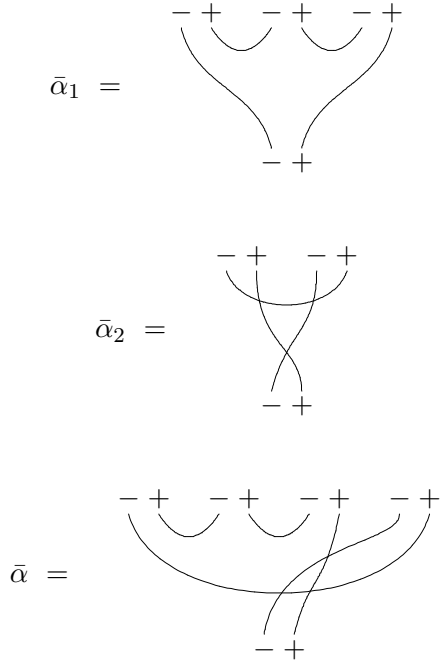

so
$$\phi_2\alpha_1 = [u \otimes u \otimes u \,,\, u]$$
$$\phi_2\alpha_2 = [u \otimes u \,,\, u]$$
$$\phi_2\alpha = [u \otimes u \otimes u \otimes u \,,\, u]$$

Then $\bar{\theta}$ may have the following graph in $A\mathbf{Ope}_1$

$$[\,u\otimes u\otimes u\,,\,u\,]\otimes[\,u\otimes u\,,\,u\,]$$

$$[\,u\otimes u\otimes u\otimes u\,,\,u\,]$$

The condition B is seen to be satisfied by the following diagram; we apply $\phi_1$ to each component, and compose with $\alpha_1 \otimes \alpha_2$:

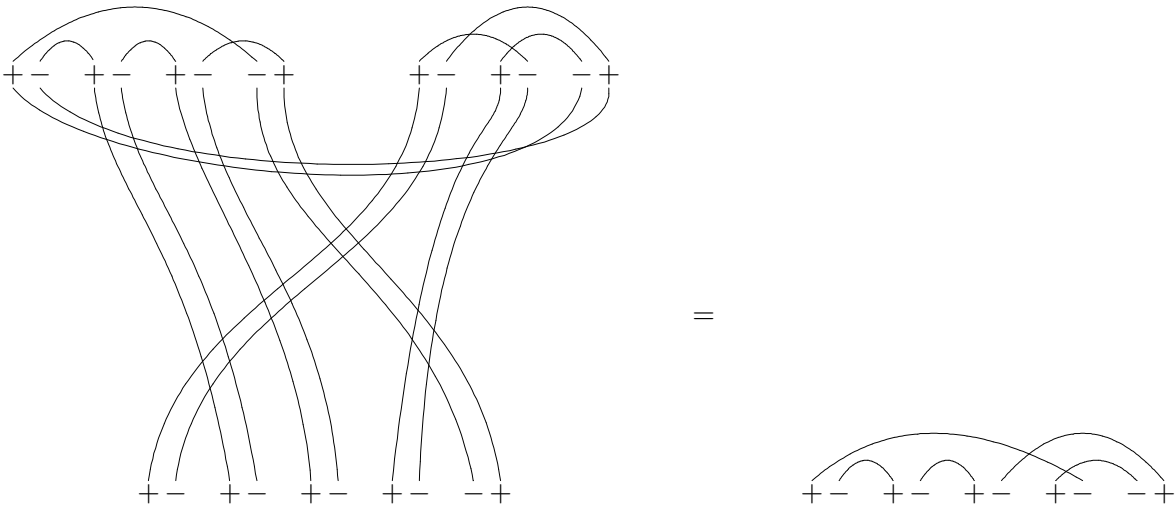

=

This corresponds to a 3-opetope of the form

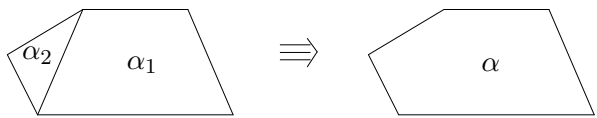

.

Note that we still do not need to label the edges of the graph since $\mathbf{Ope}_1$ also only has identity arrows.

A morphism
$$\theta \xrightarrow{f} \theta' \in \mathbf{Ope}_3$$
then has one of the following two forms

$$[\ \alpha_1 \otimes \alpha_2\ ,\ \alpha\ ]$$

$$g_1 \quad g_2 \qquad g$$

$$[\ \beta_1 \otimes \beta_2\ ,\ \beta\ ]$$

or

$$[\ \alpha_1 \otimes \alpha_2\ ,\ \alpha\ ]$$

$$g_1 \quad g_2 \qquad g$$

$$[\ \beta_1 \otimes \beta_2\ ,\ \beta\ ]$$

where $g_1, g_2, g$ are morphisms in $\mathbf{Ope}_2$. Since all morphisms in $\mathbf{Ope}_2$ are isomorphisms, it follows that all morphisms in $\mathbf{Ope}_3$ are isomorphisms. In fact, since $\mathbf{Ope}_2$ is equivalent to a discrete category, $\mathbf{Ope}_3$ is also, and similarly $\mathbf{Ope}_k$ for all $k \geq 0$; this is proved in Section 2.

- $k = 4$

Finally we give an example of a 4-opetope $\gamma \in \mathbf{Ope}_4$, with
$$\phi_4 \gamma = [\theta_1 \otimes \theta_2\ ,\ \theta]$$
where

$$\bar{\theta}_1 \ =$$

$\bar{\theta}_2 =$ 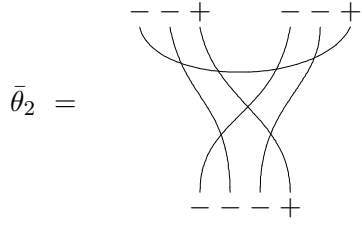

$\bar{\theta} =$ 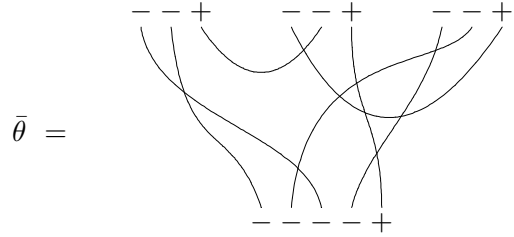

and we have

$$\phi_3\theta_1 = [\ [u^{\otimes 3}, u] \otimes [u^{\otimes 2}, u]\ ,\ [u^{\otimes 4}, u]\ ] = [U_3 \otimes U_2\ ,\ U_4], \text{ say}$$

$$\phi_3\theta_2 = [\ [u^{\otimes 2}, u] \otimes [u^{\otimes 2}, u]\ ,\ [u^{\otimes 3}, u]\ ] = [U_2 \otimes U_2\ ,\ U_3]$$

$$\phi_3\theta = [U_2 \otimes U_2 \otimes U_2\ ,\ U_4].$$

Then $\bar{\gamma}$ may be given by the following graph in $A\mathbf{Ope}_2$

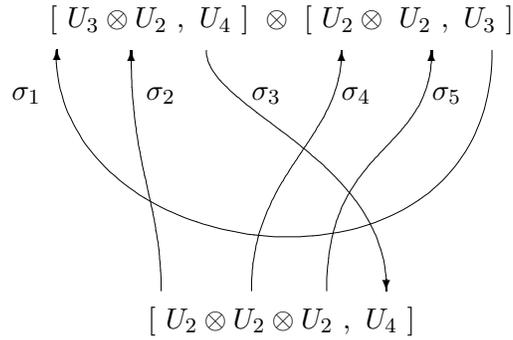

where each $\sigma_i$ is a morphism in $\mathbf{Ope}_2$, that is, a permutation. We then check condition B by the following diagram:

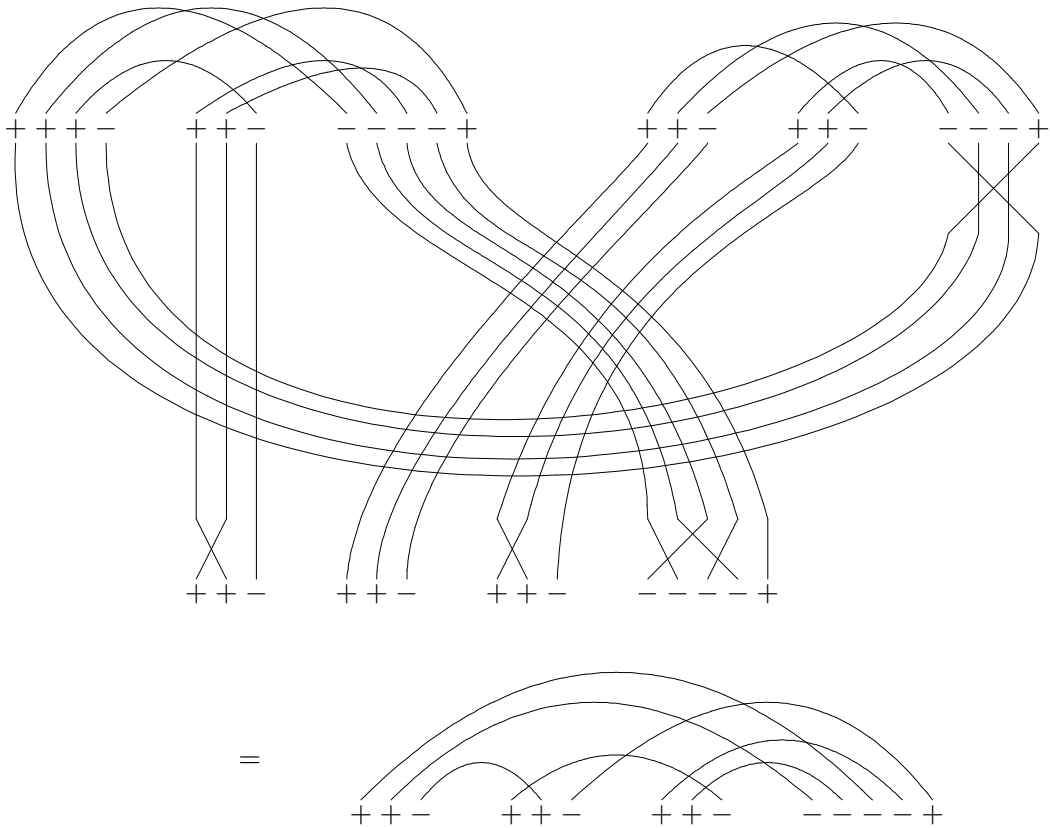

giving the composite $\theta$ as required. Note that the permutations $\sigma_i$ appear as permutations of the appropriate edges in the above diagram.

This corresponds to an opetope of the following form:

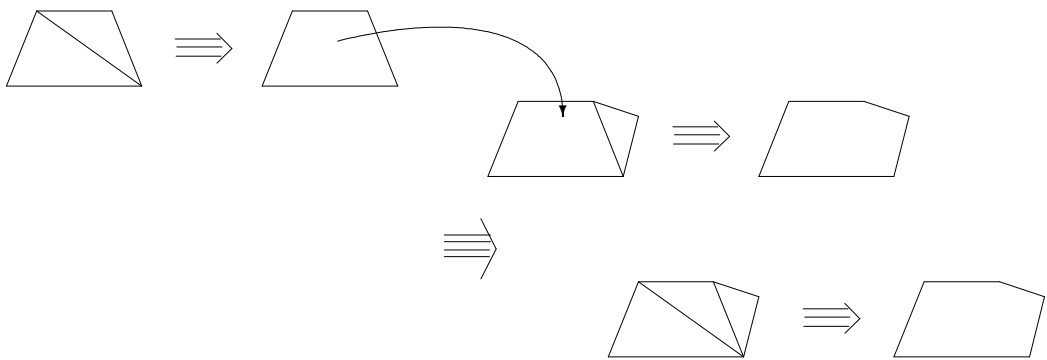

## 2 Comparison with the multicategory approach

In [1], opetopes are constructing using symmetric multicategories. Dimensions are built up using the slicing process. We compare this process with the use of closed categories as above.

### 2.1 The slice construction

Recall the slice construction for a symmetric multicategory. Let $Q$ be a symmetric multicategory. Then the slice multicategory $Q^+$ is given by

- Objects: $o(Q^+) = \text{elt}\, Q$

- Arrows: $Q^+(f_1, \ldots, f_n; f)$ is given by the set of 'configurations' for composing $f_1, \ldots, f_n$ as arrows of $Q$, to yield $f$.

Recall further that such a configuration for composing is given by a labelled tree $(T, \rho, \tau)$ where the nodes give the positions for composing the $f_i$. So by the results of [3] we may restate this using allowable morphisms in $K\mathbb{C}$, where $\mathbb{C} = o(Q)$.

Let $Q$ be a symmetric multicategory with category of objects $\mathbb{C}$. Given an arrow $f \in Q(x_1, \ldots, x_m; x)$ we write

$$\phi f = [x_1 \otimes \cdots \otimes x_m, x] \in A\mathbb{C}.$$

Then the slice multicategory $Q^+$ is given as follows.

- objects $o(Q^+) = \text{elt}\, Q$

- an arrow $\theta \in Q^+(f_1, \ldots, f_j; f)$ is an arrow

$$\theta \in A\mathbb{C}(\ I\ ,\ [\phi f_1 \otimes \cdots \otimes \phi f_j\ ,\ \phi f]\ )$$

such that composing the $f_i$ in this configuration gives $f$.

**Lemma 2.1** $\phi$ *extends to a functor*

$$\phi : \text{elt}\, Q \longrightarrow A\mathbb{C}.$$

**Proof.** Let
$$f \in Q(x_1, \ldots, x_m; x)$$
$$g \in Q(y_1, \ldots, y_j; y).$$
Then $\text{elt}\, Q(f, g) = \emptyset$ unless $m = j$. If $m = j$ then a morphism $f \xrightarrow{\gamma} g$ is given by a permutation $\sigma \in \mathbf{S}_m$ together with morphisms

$$t_i : y_i \longrightarrow x_{\sigma(i)}$$

$$t : x \longrightarrow y$$

satisfying certain conditions. This specifies a unique allowable morphism

$$[x_1 \otimes \cdots \otimes x_m, x] \longrightarrow [y_1 \otimes \cdots \otimes y_m, y] \in A\mathbb{C}$$

and we define $\phi\gamma$ to be this morphism. This makes $\phi$ into a functor. $\square$

We call $\phi$ the *frame functor* for $Q$. We now show how the slicing process corresponds to moving one rung up the 'ladder'.

**Lemma 2.2** *Let $Q$ be a symmetric multicategory with category of objects $\mathbb{C}$. Then the category $\mathrm{elt}\, Q^+$ is isomorphic to a full subcategory of the comma category $(I \downarrow A\phi)$ and the frame functor for $Q^+$ is given by*

$$\mathrm{elt}\, Q^+ \hookrightarrow (I \downarrow A\phi) \longrightarrow A(\mathrm{elt}\, Q)$$

*where the functors shown are the inclusion followed by the forgetful functor.*

**Proof.** Write $\mathbb{C}_1 = \mathrm{elt}\, Q = o(Q^+)$.

An object of $\mathrm{elt}\, Q^+$ is $(\theta, p)$ where $p \in \mathcal{F}\mathbb{C}_1^{\mathrm{op}} \times \mathbb{C}_1$ and $\theta \in Q^+(p)$.

Write

$$p = (f_1, \ldots, f_m; f).$$

Then $\theta$ is an allowable morphism

$$I \xrightarrow{\theta} A\phi[f_1 \otimes \cdots \otimes f_m, f]$$

that is, an object

$$(\,\theta\,,\, [f_1 \otimes \cdots \otimes f_m, f]\,) \in (I \downarrow A\phi)$$

such that composing the $f_i$ according to $\theta$ results in $f$.

A morphism $(\theta, p) \longrightarrow (\theta', p')$ in $\mathrm{elt}\, Q^+$ is a morphism $p \longrightarrow p'$ in $\mathcal{F}\mathbb{C}_1^{\mathrm{op}} \times \mathbb{C}_1$ such that a certain commuting condition holds. Such a morphism is precisely an allowable morphism

$$[f_1 \otimes \cdots \otimes f_m, f] \longrightarrow [f'_1 \otimes \cdots \otimes f'_m, f'] \in A\mathbb{C}_1$$

and the commuting condition is precisely that ensuring that this is a morphism $\theta \longrightarrow \theta'$ in $(I \downarrow A\phi)$.

It is then clear that the frame functor is given by the inclusion followed by the forgetful functor as asserted. $\square$

**Corollary 2.3** *The category $\mathrm{elt}\, Q^+$ is the full subcategory of $(I \downarrow A\phi)$ whose objects are all $(\theta, p)$ satisfying the following two conditions:*

i) $p$ has shape $X_m$ for some $m \geq 0$ so $p = [f_1 \otimes \cdots \otimes f_m, f]$

ii) the result of composing the $f_i$ according to $\theta$ is $f$.

If $Q$ is itself a slice multicategory then we can state the condition (ii) in the language of closed categories as well, since each $f_i$ is itself an allowable graph.

So we now consider forming $Q^{++}$, that is, the slice of a slice multicategory. Let $Q$ be a symmetric multicategory with category of objects $\mathbb{C}_0$. We write
$$\mathbb{C}_1 = o(Q^+)$$
with frame functor
$$\phi_1 : \quad \begin{array}{ccc} \mathbb{C}_1 & \longrightarrow & A\mathbb{C}_0 \\ f \in Q(x_1, \ldots, x_m; x) & \longmapsto & [x_1 \otimes \cdots \otimes x_m, x] \end{array}$$

Also, we write
$$\mathbb{C}_2 = \operatorname{elt} Q^+$$
with frame functor
$$\phi_2 : \quad \begin{array}{ccc} \mathbb{C}_2 & \longrightarrow & A\mathbb{C}_1 \\ \alpha \in Q^+(f_1, \ldots, f_m; f) & \longmapsto & [f_1 \otimes \cdots \otimes f_m, f] \end{array}$$

**Lemma 2.4** *Let $\theta$ be a configuration for composing $\alpha_1, \ldots, \alpha_j \in \operatorname{elt} Q^+ = \mathbb{C}_2$ expressed as an allowable morphism*
$$I \xrightarrow{\theta} [\, \phi_2\alpha_1 \otimes \cdots \otimes \phi_2\alpha_j \,,\, \phi_2\alpha] \in A\mathbb{C}_1.$$

*Then the result of composing the $\alpha_i$ in this configuration is*
$$(A\phi_1)\bar{\theta} \,\circ\, (\alpha_1 \otimes \cdots \otimes \alpha_j)$$

*composed as morphisms of $A\mathbb{C}_0$.*

**Proof.** By definition, each $\alpha_i$ is a morphism in $A\mathbb{C}_0$ of shape
$$I \longrightarrow [\, X_{im_1} \otimes \cdots \otimes X_{im_i} \,,\, X \,],$$
so is a tree labelled in $\mathbb{C}_0$. These trees are composed by node-replacement composition (see [3]) and the "composition graph" is given by $\bar{\theta}$. □

**Corollary 2.5** *An arrow $\theta \in Q^{++}(\alpha_1, \otimes \cdots \otimes, \alpha_j; \alpha)$ is precisely a morphism*
$$\theta \in A\mathbb{C}_1(\, I \,,\, [\phi_2\alpha_1 \otimes \cdots \otimes \phi_2\alpha_j \,,\, \phi_2\alpha] \,)$$
*such that*
$$(A\phi_1)\bar{\theta} \,\circ\, (\alpha_1 \otimes \cdots \otimes \alpha_j) = \alpha \ \in A\mathbb{C}_0$$

**Corollary 2.6** $\operatorname{elt} Q^{++}$ *is the full subcategory of* $(I \downarrow A\phi_2)$ *whose objects are all* $(\theta, p)$ *satisfying the following two conditions:*

1) $p$ has shape $X_m$ for some $m \geq 0$, so $p = [\alpha_1 \otimes \cdots \otimes \alpha_m; \alpha] \in A\mathbb{C}_2$

2) $(A\phi_1)\bar{\theta} \circ (\alpha_1 \otimes \cdots \otimes \alpha_j) = \alpha$.

Finally we are ready to show that the opetopes constructed using symmetric multicategories correspond to those constructed in closed categories.

**Corollary 2.7** *Let $Q$ be the symmetric multicategory with just one object and one (identity) morphism. Then for all $k \geq 0$*

$$o(Q^{k+}) \cong \mathbf{Ope}_k$$

*where* $Q^{0+} = Q$.

**Proof.** For $k \leq 1$ the result is immediate by Definition 1.3. For $k = 2$ we use Corollary 2.3 on $Q^+$; the result follows since condition (ii) is trivially satisfied. For $k \geq 3$ we use Corollary 2.6 on $Q^{(k-3)+}$; the result follow since the $\phi_2$ in the Corollary is $\phi_{k-2}$ in the case in question. □

## 3   The category of opetopes

Recall that in [2] we defined the category $\mathcal{O}$ of opetopes. It is now possible to restate this definition in the framework of Kelly-Mac Lane graphs; we copy the definition exactly, using the fact that the bijection giving the formal definition of a tree gives the mates in the corresponding Kelly-Mac Lane graph.

Although we do not give the construction explicitly here, we give some examples of low-dimensional face maps. We use the example of a 3-opetope as given in Section 1.3.

For the 2-opetopes we have face maps

$$\begin{array}{c}
\alpha_1 \\
\phantom{x} \searrow^{s_1} \\
\alpha_2 \xrightarrow{s_2} \theta \\
\phantom{x} \nearrow_{t} \\
\alpha
\end{array}$$

together with the isomorphic cases.

For 1-opetopes we then have

$$s_1, s_2, s_3, t \; : \; u \longrightarrow \alpha_1$$

$$s_1, s_2, t \; : \; u \longrightarrow \alpha_2$$

$$s_1, s_2, s_3, s_4, t \; : \; u \longrightarrow \alpha$$

but by considering the generating relations, here given by mates in the graph $\theta$, we have

$$\begin{aligned} s_1 s_1 &= s_2 t \\ s_1 s_2 &= t s_2 \\ s_1 s_3 &= t s_3 \\ s_1 t &= tt \\ s_2 s_1 &= t s_1 \\ s_2 s_2 &= t s_4; \end{aligned}$$

note that $s_i s_j$ give the $jth$ source of the $ith$ source of $\theta$.

For 0-opetopes we have in addition face maps

$$x \longrightarrow u$$

and the relations on composites

$$x \longrightarrow \theta$$

are generated by relations on composites

$$x \longrightarrow \alpha_i$$

as well as by those on composites

$$u \longrightarrow \theta.$$

For the former relations we are considering mates under graphs $\alpha_i \in A\mathbf{Ope}_0$, and for the latter, mates under the graph $(A\phi_1)\bar{\theta} \in A\mathbf{Ope}_0$. So in fact we are considering, in total, all objects connected in the composite graph

$$(A\phi_1)\bar{\theta} \circ (\alpha \otimes \cdots \otimes \alpha_m) \in A\mathbf{Ope}_1.$$

So we have

$$\begin{aligned} ts_1 s &= s_2 s_1 s = s_2 s_2 t = ts_4 t \\ ts_1 t &= s_2 s_1 t = s_2 tt = s_1 s_1 t = s_1 s_2 s = ts_2 s \\ ts_2 t &= s_1 s_2 t = s_1 s_3 s = ts_3 s \\ ts_3 t &= s_1 s_3 t = s_1 tt = ttt \\ ts_4 s &= s_2 s_2 s = s_2 ts = s_1 s_1 s = s_1 ts = tts \end{aligned}$$

Note that since

$$(A\phi_1)\bar{\phi} \circ (\alpha_1 \otimes \cdots \otimes \alpha_m) = \alpha$$

the 0-cell face maps for $\theta$ are precisely those of the form $tf$ where $f$ is a 0-cell face map of $\alpha = t(\theta)$. This reflects the fact that, when 2-opetopes are composed along 1-opetopes, the composite is formed by 'deleting' the boundary 1-opetopes, but no 0-cells are deleted. This result generalises to $k$-opetopes, but we do not prove this here.